\documentclass[11pt]{article}%
\usepackage{amsfonts}
\usepackage{amsmath}
\usepackage{amssymb}
\usepackage{amscd}
\usepackage{graphicx}
\usepackage{amsthm}%
\newtheorem{definition}{Definition}[section]
\newtheorem{proposition}{Proposition}[section]
\newtheorem{theorem}{Theorem}[section]
\theoremstyle{definition}
\newtheorem{example}{Example}[section]
\newenvironment{remarks}{{\it Remarks.}}{\hfill $\triangleleft$ \smallskip}

\renewcommand{\[}{\begin{equation}}
\renewcommand{\]}{\end{equation}}

\begin{document}

\title{Kernel Methods for Linear Discrete-Time Equations}
\date{}
\author{Fritz Colonius\\Institut f\"{u}r Mathematik, Universit\"{a}t Augsburg, Augsburg/Germany
\and Boumediene Hamzi\\Department of Mathematics, Ko\c{c} University, Istanbul, Turkey.}
\maketitle

\textbf{Abstract}: Methods from learning theory are used in the state space of
linear dynamical and control systems in order to estimate the system matrices.
An application to stabilization via algebraic Riccati equations is included.
The approach is illustrated via a series of numerical examples.


\section{Introduction}

This paper discusses several problems in dynamical systems and control, where
methods from learning theory are used in the state space of linear systems.
This is in contrast to previous approaches in the frequency domain \cite{zhou,
ljung}. We refer to \cite{ljung} for a general survey on applications of
machine learning to system identification.

Basically, learning theory allows to deal with problems when only data from a
given system are given. Reproducing Kernel Hilbert Spaces (RKHS) allow to work
in a very large dimensional space in order to simplify the underlying problem.
We will discuss this in the simple case when the matrix $A$ describing a
linear discrete-time system is unknown, but a time series from the underlying
linear dynamical system is given. We propose a method to estimate the
underlying matrix using kernel methods. Applications are given in the stable
and unstable case and for estimating the topological entropy for a linear map.
Furthermore, in the control case, stabilization via linear-quadratic optimal
control is discussed.

The emphasis of the present paper is on the formulation of a number of
problems in dynamical systems and control and to illustrate the applicability
of our approach via a series of numerical examples.

The contents is as follows: In Section \ref{Section2} the problem is stated
formally and an algorithm based on kernel methods is given for the stable
case. In Section \ref{Section3} the algorithm is extended to the unstable
case. In particular, the topological entropy of linear maps is computed (which
boils down to computing unstable eigenvalues). In Section \ref{Section4b}
identification of linear control systems is considered and Section
\ref{Section5} discusses their stabilization. Here we insert the estimate of
the system matrix (obtained via learning theory) into the relevant algebraic
Riccati equation and study when this yields a stabilizing feedback. Every
section contains several numerically computed examples (via MATLAB)
illustrating the approach. Section \ref{Conclusions} draws some conclusions
from the numerical experiments. For the reader's convenience we have collected
in the appendix basic concepts from learning theory as well as some hints to
the relevant literature.

\section{Statement of the problem\label{Section2}}

Consider the linear discrete-time system
\[
\label{syst1}x(k+1)=Ax(k),
\]
where $A=[a_{i,j}]\in\mathbb{R}^{n\times n}$. We want to estimate $A$ from the
time series $x(1)+\eta_{1}$, $\cdots$, $x(N)+\eta_{N}$ where the initial
condition $x(0)$ is known and $\eta_{i}$ are distributed according to a
probability measure $\rho_{x}$ that satisfies the following condition (this is
the \emph{Special Assumption} in \cite{smale_shannon1}).

\textbf{Assumption } The measure $\rho_{x}$ is the marginal on $X=\mathbb{R}%
^{n}$ of a Borel measure $\rho$ on $X\times\mathbb{R}$ with zero mean
supported on $[-M_{x},M_{x}],M_{x}>0$.

One obtains from (\ref{syst1}) for the components of the time series that
\[
\label{param1}x_{i}(k+1)=\sum_{j=1}^{n}a_{ij}x_{j}(k).
\]
For every $i$ we want to estimate the coefficients $a_{ij},j=1,$ $\cdots,n$.
They are determined by the linear maps $f^{\ast}_{i}:\mathbb{R}^{n}%
\rightarrow\mathbb{R}$ given by%
\[
\label{unknown_map}(x_{1},...,x_{n})\mapsto\sum_{j=1}^{n}a_{ij}x_{j}.
\]
This problem can be reformulated as a learning problem as described in the
Appendix where $f^{\ast}_{i}$ in (\ref{unknown_map}) plays the role of the
unknown function (\ref{unknown_func}) and $(x(k),x_{i}(k+1)+\eta_{i})$ are the
samples in (\ref{samples}).

We note that in \cite{smale_shannon1}, the authors do not consider time series
and that we apply their results to time series.

In order to approximate $f^{\ast}_{i}$, we minimize the criterion in
(\ref{reg_opt2}). For a positive definite kernel $K$ let $f_{i}$ be the kernel
expansion of $f^{\ast}_{i}$ in the corresponding RKHS $\mathcal{H}_{K}$. Then
$f_{i}=\sum_{j=1}^{\infty}c_{i,j}\phi_{j}$ with certain coefficients
$c_{ij}\in\mathbb{R}$ and%
\[
||f_{i}||_{\mathcal{H}_{K}}=\displaystyle\sum_{j=1}^{\infty}\frac{c_{i,j}^{2}%
}{\lambda_{j}},
\]
where $(\lambda_{j},\phi_{j})$ are the eigenvalues and eigenfunctions of the
integral operator $L_{K}:\mathcal{L}_{\nu}^{2}(\mathcal{X})\rightarrow
\mathcal{C}(\mathcal{X})$ given by $(L_{K}f)(x)=\int K(x,t)f(t)d\nu(t)$ with a
Borel measure $\nu$ on $\mathcal{X}$. Thus $L_{K}\phi_{j}=\lambda_{j}\phi_{j}$
for $j\in\mathbb{N}^{\ast}$ and the eigenvalues $\lambda_{j}\geq0$.

Then we consider the problem of minimizing over $(c_{i,1},$ $\cdots,c_{i,N})$
the functional%

\[
\label{error}\mathcal{E}_{i}=\frac{1}{N}\sum_{k=1}^{N}(y_{i}(k)-f_{i}%
(x(k)))^{2}+\gamma_{i}||f_{i}||_{\mathcal{H}_{K}}^{2},
\]
where $y_{i}(k):=x_{i}(k+1)+\eta_{i}=f^{\ast}_{i}(x(k))+\eta_{i}$ and
$\gamma_{i}$ is a regularization parameter.


Since we are dealing with a linear problem, it is natural to choose the linear
kernel $k(x,y)=\langle x,y\rangle$. Then the solution of the above
optimization problem is given by the kernel expansion of $x_{i}(k+1)$,
$i=1,\cdots,n$,%

\[
\label{cij}y_{i}(k):=x_{i}(k+1)=\sum_{j=1}^{N}c_{ij}\langle x(j),x(k)\rangle,
\]
where the $c_{ij}$ satisfy the following set of equations:%

\begin{equation}
\left[
\begin{array}
[c]{c}%
x_{i}(1)\\
\vdots\\
x_{i}(N)
\end{array}
\right]  =\Bigg(N\lambda I_{d}+{\mathbb{K}}\Bigg)\left[
\begin{array}
[c]{c}%
c_{i1}\\
\vdots\\
c_{iN}%
\end{array}
\right]  , \label{cij_eqn}%
\end{equation}
with
\[
{\mathbb{K:}}=\left[
\begin{array}
[c]{ccc}%
\sum_{\ell=1}^{n}x_{\ell}(1)x_{\ell}(0) & \cdots & \sum_{\ell=1}^{n}x_{\ell
}(N)x_{\ell}(0)\\
\vdots & \cdots & \vdots\\
\sum_{\ell=1}^{n}x_{\ell}(1)x_{\ell}(N-1) & \cdots & \sum_{\ell=1}^{n}x_{\ell
}(N)x_{\ell}(N-1)
\end{array}
\right]  .
\]
This is a consequence of Theorem \ref{thm_sol1}.

From (\ref{param1}), we have
\begin{align*}
x_{i}(k+1)  &  =\sum_{j=1}^{N}c_{ij}\langle x(j),x(k)\rangle=\sum_{j=1}%
^{N}c_{ij}x(j)^{T}\cdot x(k)=\sum_{j=1}^{N}\sum_{\ell=1}^{n}c_{ij}x_{\ell
}(j)x_{\ell}(k)\\
&  =\sum_{\ell=1}^{n}\sum_{j=1}^{N}c_{ij}x_{\ell}(j)x_{\ell}(k).
\end{align*}
Then an estimate of the entries of $A$ is given by
\begin{equation}
\hat{a}_{i\ell}=\sum_{j=1}^{N}c_{i,j}x_{\ell}(j). \label{aij_eqn}%
\end{equation}
This discussion leads us to the following basic algorithm.

\paragraph{Algorithm $\mathcal{A} $:}

If the eigenvalues of $A$ are all within the unit circle, one proceeds as
follows in order to estimate $A$. Given the time series $x(1),\cdots,x(N)$
solve the system of equations (\ref{cij_eqn}) to find the numbers $c_{ij}$ and
then compute $\hat{a}_{i\ell}$ from (\ref{aij_eqn}).\medskip

Before we present numerical examples and modifications and applications of
this algorithm, it is worthwhile to note the following preliminary remarks
indicating what may be expected.

The stability assumption in algorithm $\mathcal{A}$ is imposed, since
otherwise the time series will diverge exponentially. Then, already for a
moderately sized number of data points ($N\approx10^{2}$) equation
(\ref{cij_eqn}) will be ill conditioned. Hence for unstable $A$, modifications
of algorithm $\mathcal{A}$ are required.

While for test examples one can compare the entries of the matrix $A$ and its
approximation $\hat{A}$, it may appear more realistic to compare the values
$x(1),\cdots,x(N)$ of the data series and the values $\hat{x}(1),\cdots
,\hat{x}(N)$ generated by the iteration of the matrix $\hat{A}$.

In general, one should not expect that increasing the number of data points
will lead to better approximations of the matrix $A$. If the matrix $A$ is
diagonalizable, for generic initial points $x(0)\in\mathbb{R}^{n}$ the data
points $x(k)$ will approach for $N\rightarrow\infty$ the eigenspace for the
eigenvalue with maximal modulus. For general $A$ and generic initial points
$x(0)\in\mathbb{R}^{n}$, the data points $x(N)$ will approach for
$N\rightarrow\infty$ the largest Lyapunov space (i.e., the sum of the real
generalized eigenspaces for eigenvalues with maximal modulus). Thus in the
limit for $N\rightarrow\infty$, only part of the matrix can be approximated. A
detailed discussion of this (well known) limit behavior is, e.g., given in
Colonius and Kliemann \cite{ColK00}. A consequence is that a medium length of
the time series should be adequate.

This problem can be overcome by choosing the regularization parameter $\gamma$
in (\ref{error}) and (\ref{cij_eqn}) using the method of cross validation
described in \cite{rifkin}. Briefly, in order to choose $\gamma$, we consider
a set of values of regularization parameters: we run the learning algorithm
over a subset of the samples for each value of the regularization parameter
and choose the one that performs the best on the remaining data set. Cross
validation helps also in the presence of noise and to improve the results
beyond the training set.

A theoretical justification of our algorithm could be guaranteed by the error
estimates in Theorem \ref{thm:errors}. In fact, for the linear dynamical
system (\ref{syst1}), we have that $f^{\ast}$ in (\ref{unknown_func}) is the
linear map $f^{\ast}(x)=f_{i}(x)$ in (\ref{unknown_map}) and the samples
$\mathbf{s}$ in (\ref{samples}) are $(x(k),x_{i}(k+1)+\eta_{i})$. Moreover, by
choosing the linear kernel $k(x,y)=\langle x,y\rangle$ we get that $f^{\ast
}\in\mathcal{H}_{K}$. 


Next we discuss several numerical examples, beginning with the following
scalar equation.

\begin{example}
\label{Example1}Consider $x(k+1)=\alpha x(k)$ with $\alpha=0.5$. With the
initial condition $x(0)=-0.5$, we generate the time series $x(1),\cdots
,x(100)$. Applying algorithm $\mathcal{A}$ with the regularization parameter
$\gamma=10^{-6}$ we compute $\hat{\alpha}=0.4997$. Using cross validation, we
get that $\hat{\alpha}=0.5$ with regularization parameter $\gamma
=1.5259\cdot10^{-5}$. When we introduce an i.i.d perturbation signal $\eta
_{i}\in\lbrack-0.1,0.1]$, the algorithm does not behave well when we fix the
regularization parameter. With cross validation, the algorithm works quite
well and the regularization parameter adapts to the realization of the signal
$\eta_{i}$. Here, for $e(k)=x(k)-\hat{x}(k)$ with $x(k+1)=\alpha x(k)$ and
$\hat{x}(k+1)=\hat{\alpha}\hat{x}(k)$, we get that $||e(300)||=\sqrt
{\sum_{i=1}^{300}e^{2}(i)}=0.0914$ and $\sqrt{\sum_{i=100}^{300}e^{2}%
(i)}=1.8218\cdot10^{-30}$.

We observe an analogous behavior of the algorithm when the data are generated
from $x(k+1)=\alpha x(k)+\varepsilon x(k)^{2}$ where the algorithm works well
in the presence of noise and structural perturbations when using cross
validation. When $\varepsilon=0.1$ and with an i.i.d perturbation signal
$\eta_{i}\in\lbrack-0.1,0.1]$, $\hat{\alpha}$ varies between $0.38$ and $0.58$
depending on the realization of $\eta_{i}$ but $||e(300)||=\sqrt{\sum
_{i=1}^{300}e^{2}(i)}=0.2290$ and $\sqrt{\sum_{i=100}^{300}e^{2}%
(i)}=2.8098\cdot10^{-30}$ which shows that the error $e$ decreases
exponentially and the generalization properties of the algorithm are quite good.
\end{example}

\begin{example}
\label{Example3}Consider $x(k+1)=Ax(k)$ with matrix $A$ given by
\[
A:=\left[
\begin{array}
[c]{cccc}%
-0.5 & 1 & 0 & 0\\
0 & 0.6 & 1 & 0\\
0 & 0 & 0.7 & 1\\
0 & 0 & 0 & -0.8
\end{array}
\right]  .
\]
For the initial condition $x=[-0.9,0.1,15,0.2]^{\prime}$ and with $N=100$ data
points, we get%
\[
\hat{A}=\left[
\begin{array}
[c]{cccc}%
-0.5000 & 1.0000 & 0.0000 & -0.0000\\
0.0000 & 0.6000 & 1.0000 & 0.0000\\
0.0000 & -0.0000 & 0.7000 & 0.9994\\
-0.0000 & 0.0000 & -0.0000 & -0.7995
\end{array}
\right]  .
\]
We then simulate $x(k+1)=Ax(k)$ and $\hat{x}(k+1)=\hat{A}\hat{x}(k)$ for
$k=0,\cdots,200$ to test the accuracy of our approximation beyond the interval
$k=0,\cdots,100$. Then the norm of the error $e_{j}(k)=x_{j}(k)-\hat{x}%
_{j}(k)$, for $j=1,\cdots,4$, $||e_{j}(300)||=\sqrt{\sum_{i=1}^{300}e_{j}%
^{2}(i)}$ is of the order of $10^{-3}$ and $\sqrt{\sum_{i=100}^{300}e_{j}%
^{2}(i)}$ is of the order of $10^{-11}$ which shows that the error $e$
decreases exponentially and the generalization properties of the algorithm are
quite good. The regularization parameters are $\gamma_{i}=0.9313\cdot10^{-9}$
for $i=1,\cdots,4$.

Also in the presence of small noise $\eta_{i}\in\lbrack-0.01,0.01]$, the
algorithm behaves well and the regularization parameters adapt to the
realization of $\eta_{i}$. For example, for a certain realizations of
$\eta_{i}$, we obtain the regularization parameters%
\[
\gamma_{1}=0.0039, \gamma_{2}=2.4114\cdot10^{-4}, \gamma_{3}=9.3132\cdot
10^{-10}, \gamma_{4}=2 \cdot10^{-3}
\]
and the error $||e_{j}(300)||=\sqrt{\sum_{i=1}^{300}e_{j}^{2}(i)}$ is of the
order of $10^{-1}$ and $\sqrt{\sum_{i=100}^{300}e_{j}^{2}(i)}$ is of the order
of $10^{-9}$ .

Suppose that in addition to a small noise $\eta_{i}\in\lbrack-0.01,0.01],$
there is a quadratic structural perturbation, i.e.,%
\[
x(k+1)=Ax(k)+\varepsilon\left[
\begin{array}
[c]{c}%
x_{1}(k)^{2}\\
x_{2}(k)^{2}\\
x_{3}(k)^{2}\\
x_{4}(k)^{2}%
\end{array}
\right]  .
\]
Then with cross validation for $\varepsilon=0.001$ the algorithm behaves well.
For a particular realization of $\eta$, the error $||e_{j}(300)||=\sqrt
{\sum_{i=1}^{300}e_{j}^{2}(i)}$ is between $5$ and $15$ but $\sqrt
{\sum_{i=100}^{300}e_{j}^{2}(i)}$ is of the order of $10^{-9}$ and the
regularization parameters are
\[
\gamma_{1}=0.5,\gamma_{2}=9.3132\cdot10^{-10},\gamma_{3}=9.3132\cdot
10^{-10},\gamma_{4}=9.3132\cdot10^{-10}.
\]

\end{example}

These examples show a very good behavior of the algorithm.

\section{Unstable case\label{Section3}}

Consider
\begin{equation}
x(k+1)=Ax(k)\text{ with }A\in\mathbb{R}^{n\times n}, \label{0}%
\end{equation}
where some of the eigenvalues of $A$ are outside the unit circle. Again, we
want to estimate $A$ when the following data are given,
\begin{equation}
x(1),x(2),...,x(N), \label{1}%
\end{equation}
which are generated by system (\ref{0}), thus $x(k)=A^{k-1}x(1)$.

As remarked above, a direct application of the algorithm $\mathcal{A}$ will
not work, since the time series diverges fast. Instead we construct a new time
series from (\ref{1}) associated to an auxiliary stable system.

For a constant $\sigma>0$ we define the auxiliary system by
\begin{equation}
y(k+1)=\tilde{A}y(k)\text{ with }\tilde{A}:=\frac{1}{\sigma}A. \label{2}%
\end{equation}
Thus%
\[
y(k)=\left(  \frac{A}{\sigma}\right)  ^{k-1}y(1)
\]
and with $y(1)=x(1)$ one finds%
\[
y(k)=\frac{1}{\sigma^{k-1}}A^{k-1}x(1)=\frac{1}{\sigma^{k-1}}x(k).
\]
If we choose $\sigma>0$ such that the eigenvalues of $\frac{A}{\sigma}$ are in
the unit circle, we can apply algorithm $\mathcal{A}$ to this stable matrix
and hence we would obtain an estimate of $\frac{A}{\sigma}$ and hence of $A$.
However, since the eigenvalues of the matrix $A$ are unknown, we will be
content with a somewhat weaker condition than stability of $\frac{A}{\sigma}$.

The data (\ref{1}) for system (\ref{0}) yield the following data for system
(\ref{2}):%
\[
y(1):=x(1),y(2):=\frac{1}{\sigma}x(2),...,y(N):=\frac{1}{\sigma^{N-1}}x(N).
\]
We propose to choose $\sigma$ as follows: Define
\[
\label{gamma}\sigma:=\max\left\{  \frac{\left\Vert x(k+1)\right\Vert
}{\left\Vert x(k)\right\Vert },k\in\{0,1,...,N\}\right\}  .
\]
Clearly the inequality $\sigma\leq\left\Vert A\right\Vert $ holds. We apply
algorithm $\mathcal{A}$ to the time series $y(k)$. This yields an estimate of
$\frac{A}{\sigma}$ and hence an estimate $\hat{A}$ of $A$.

For general $A$, this choice of $\sigma$ certainly does not guarantee that the
eigenvalues of $\frac{A}{\sigma}$ are within the unit circle. However, as
mentioned above, a generic data sequence $x(k),k\in\mathbb{N}$, will converge
to the eigenspace of the eigenvalue with maximal modulus. Hence $\frac
{\left\Vert x(k+1)\right\Vert }{\left\Vert x(k)\right\Vert }$ will approach
the maximal modulus of an eigenvalue, thus this choice of $\sigma$ will lead
to a matrix $\frac{A}{\sigma}$ which is not \textquotedblleft too
unstable\textquotedblright.

\begin{example}
\label{Example4}Consider $x(k+1)=\alpha x(k)$ with $\alpha=11.46$. With the
initial condition $x(0)=-0.5$, we generate the time series $x(1),\cdots
,x(100)$. The algorithm above with the regularization parameter $\gamma
=10^{-6}$ yields the estimate $\hat{\alpha}=11.4086$. Cross validation leads
to the regularization parameter $\gamma=9.5367\cdot10^{-7}$ and the estimate
$\hat{\alpha}=11.4599$. In the presence of a small noise $\eta\in
\lbrack-0.1,0.1]$, cross validation yields the regularization parameter
$\gamma=0.002$ and the slightly worse estimate $\hat{\alpha}=11.1319$.
\end{example}

We observe the same behavior in higher dimensional systems where the
eigenvalues are of the same order of magnitude.

\begin{example}
\label{Example4b}Consider $x(k+1)=Ax(k)$ with%
\[
A=\left[
\begin{array}
[c]{cccc}%
20 & 0 & 0 & 0\\
0 & -10 & 0 & 0\\
0 & 0 & 15 & 0\\
0 & 0 & 0 & -25
\end{array}
\right]
\]
Using cross validation, we get that%
\[
\hat{A}=\left[
\begin{array}
[c]{cccc}%
20.0000 & 0.0000 & 0.0001 & 0.0000\\
-0.0000 & -10.0000 & 0.0000 & -0.0000\\
0.0000 & -0.0000 & 14.9998 & 0.0000\\
-0.0000 & -0.0000 & -0.0000 & -25.0003
\end{array}
\right]
\]
for $\gamma_{i}=0.9313\cdot10^{-9}$, $i=1,\cdots,4$.

For different realizations of a noise $\eta_{i}$ of magnitude $0.5\cdot
10^{-4}$, cross validation gives a good approximation of $A$ and the
eigenvalues of $A-\hat{A}$ are all within the unit disk with amplitude of the
order of $10^{-3}$ showing that the dynamics of the error $e(k)=x(k)-\hat
{x}(k)$ is asymptotically stable. For example, for a particular realization of
$\eta_{i}$ of magnitude $0.5\cdot10^{-4}$, we get
\[
\hat{A}=\left[
\begin{array}
[c]{cccc}%
19.9635 & 0.0086 & 0.1365 & -0.0007\\
-0.0177 & -10.0025 & 0.0379 & -0.0007\\
-0.0177 & -0.0025 & 15.0376 & -0.0007\\
-0.0132 & -0.0167 & 0.0065 & -25.0000
\end{array}
\right]
\]
with regularization parameters%
\[
\gamma_{1}=1.9073\cdot10^{-6},\gamma_{2}=9.3132\cdot10^{-10},\gamma
_{3}=9.3132\cdot10^{-10},\gamma_{4}=1.2207\cdot10^{-4}.
\]
The algorithm fails in the presence of quadratic structural perturbations.
This is due to the choice of a linear kernel. A polynomial kernel, for
example, would allow for nonlinear perturbations but this would require a
complete reformulation of our algorithm. We leave the extension of our
algorithm to the nonlinear case for future work.
\end{example}

The next example is an unstable system with a large gap between the eigenvalues.

\begin{example}
\label{Example5}Consider the system $x(k+1)=Ax(k)$ with
\[
A=\left[
\begin{array}
[c]{cc}%
20 & 0\\
0 & -0.1
\end{array}
\right]  .
\]
With the initial condition $x(0)=[-1.9,1]$, we generate the time series
$x(1),\cdots,x(100)$. The algorithm above yields the (excellent) estimate
\[
\hat{A}=\left[
\begin{array}
[c]{cc}%
20.0000 & 0.0000\\
-0.0000 & -0.1000
\end{array}
\right]  ,
\]
In the presence of noise of maximal amplitude $10^{-4}$ , the algorithm
approximates well only the large entry $a_{11}=20$: For a first realization of
$\eta_{i}$ and with cross validation, we get
\[
\hat{A}=\left[
\begin{array}
[c]{cc}%
19.9997 & -0.0111\\
0.0000 & -0.1104
\end{array}
\right]  ,
\]
with $\gamma_{1}=1.5259\cdot10^{-5}$ and $\gamma_{2}=2^{20}$. However another
realization of $\eta_{i}$ leads to
\[
\hat{A}=\left[
\begin{array}
[c]{cc}%
19.9994 & -0.0011\\
0.0000 & -0.0000
\end{array}
\right]  ,
\]
with $\gamma_{1}=3.0518\cdot10^{-5}$ and $\gamma_{2}=2.8147\cdot10^{14}$. This
is due to the fact that the data converge to the eigenspace generated by the
largest eigenvalue $\lambda=20$. However, the eigenvalues of $A-\hat{A}$ are
within the unit disk with small amplitude which guarantees that the error
dynamics of $e(k)=x(k)-\hat{x}(k)$ converges to the origin quite quickly. We
observe the same phenomenon with
\[
A=\left[
\begin{array}
[c]{cc}%
-0.5 & 0\\
0 & 25
\end{array}
\right]  .
\]
Here, in the absence of noise, we obtain the estimate
\[
\hat{A}=\left[
\begin{array}
[c]{cc}%
-0.5000 & 0.0000\\
-0.0000 & 25.0000
\end{array}
\right]  ,
\]
with $\gamma_{1}=\gamma_{2}=0.9313\cdot10^{-9}$. In the presence of noise
$\eta_{i}$ with amplitude $10^{-4}$, the data converge to the eigenspace
corresponding to the largest eigenvalue $\lambda=25$: for some realization of
$\eta_{i}$ one obtains the estimate%
\[
\hat{A}=\left[
\begin{array}
[c]{cc}%
-0.4809 & 0.0008\\
0.0164 & 24.9960
\end{array}
\right]  ,
\]
while for another realization of $\eta$
\[
\hat{A}=\left[
\begin{array}
[c]{cc}%
-0.0000 & -0.0000\\
-1.0067 & 24.8696
\end{array}
\right]  .
\]
The regularization parameters $\gamma_{1}$ and $\gamma_{2}$ adapt to the
realization of the noise.
\end{example}

As already remarked in the end of Section \ref{Section2}, we see that
\textquotedblleft more data\textquotedblright\ does not always necessarily
lead to better results, since the data sequence converges to the eigenspace
generated by the largest eigenvalue. However, whether with or without noise,
the approximations of ${A}$ are good enough to reduce the error between
$x(k+1)=Ax(k)$ and $\hat{x}(k+1)=\hat{A}\hat{x}(k)$ outside of the training
examples, since cross-validation determines a good regularization parameter
$\gamma$ that balances between good fitting and good prediction properties.

The next example has an eigenvalue on the unit circle.

\begin{example}
\label{Example6}Consider $x(k+1)=Ax(k)$ with
\[
A=\left[
\begin{array}
[c]{cccc}%
2.2500 & -1.2500 & 1.2500 & -49.5500\\
3.7500 & -2.7500 & 13.1500 & -20.6500\\
0 & 0 & 10.4000 & -32.3000\\
0 & 0 & 0 & -21.9000
\end{array}
\right]  .
\]
The set of eigenvalues of $A$ is
$\mbox{spec}(A)=\{-1.5000,1.0000,10.4000,-21.9000\}$. In the absence of noise
and initial condition $x=[-0.9,15,1.5.2.5]$ with $N=100$ points, we compute
the estimate%
\[
\hat{A}=\left[
\begin{array}
[c]{cccc}%
2.2500 & -1.2500 & 1.2498 & -49.5499\\
3.7500 & -2.7500 & 13.1498 & -20.6499\\
0.0000 & 0.0000 & 10.3998 & -32.2999\\
0.0000 & 0.0000 & -0.0001 & -21.8999
\end{array}
\right]  ,
\]
and regularization parameters $\gamma_{1}=\gamma_{2}=0.9313\cdot10^{-9}$. In
this case, the set of eigenvalues of $\hat{A}$ is
\[
\mbox{spec}(\hat{A})=\{-21.9000,10.3999,-1.5000,1.0000\}.
\]
For a given realization of $\eta\in\lbrack-10^{-4},10^{-4}]$, we obtain the
estimate
\[
\hat{A}=\left[
\begin{array}
[c]{cccc}%
2.2551 & -1.2490 & 1.2187 & -49.5304\\
3.7554 & -2.7489 & 13.1175 & -20.6297\\
0.0055 & 0.0011 & 10.3669 & -32.2794\\
0.0053 & 0.0010 & -0.0325 & -21.8797
\end{array}
\right]
\]
with $\gamma_{1}=0.0745\cdot10^{-7}$ and $\gamma_{2}=0.1490\cdot10^{-7}$. The
eigenvalues of $A-\hat{A}$ are of the order of $10^{-4}$ which guarantees that
the error dynamics converges quickly to the origin. However, the set of
eigenvalues of $\hat{A}$ is
\[
\mbox{spec}(\hat{A})=\{-21.8996,10.3999,-1.5026,1.0134\}.
\]
Hence an additional unstable eigenvalue occurs.
\end{example}



\begin{example}
\label{Example7}Consider $x(k+1)=Ax(k)$ with
\[
A=\left[
\begin{array}
[c]{cccc}%
-0.8500 & 0.4500 & -0.4500 & -77.8500\\
-1.3500 & 0.9500 & 14.3500 & -11.6500\\
0 & 0 & 15.3000 & -55.3000\\
0 & 0 & 0 & -40.0000
\end{array}
\right]  .
\]
The eigenvalues of $A$ are given by
\[
\mbox{spec}(A)=\{-0.4000,0.5000,15.3000,-40.0000\}.
\]
For an initial condition $x=[-0.9;15;1.5;2.5]$ and with $N=100$ data points,
we get%
\[
\hat{A}=\left[
\begin{array}
[c]{cccc}%
-0.8498 & 0.4501 & -0.4499 & -77.8504\\
-1.3499 & 0.9500 & 14.3501 & -11.6502\\
0.0001 & 0.0001 & 15.3001 & -55.3004\\
-0.0004 & -0.0002 & -0.0004 & -39.9987
\end{array}
\right]
\]
with eigenvalues given by
\[
\mbox{spec}(\hat{A})=\{-40.0000,-0.3974,0.4982,15.3008\}.
\]
Here we used $\gamma_{i}=10^{-12}$, $i=1,\cdots,4$. Moreover, the eigenvalues
of $A-\hat{A}$ are quite small and such that the error dynamics converges
quickly to the origin. In the presence of noise $\eta$, the algorithm
approximates the largest eigenvalues of $A$ but does not approximate the
smaller (stable) ones. For example, for a particular realization of noise with
amplitude $10^{-4}$, we get the estimate%
\[
\hat{A}=\left[
\begin{array}
[c]{cccc}%
-2.1100 & -0.0993 & -1.3259 & -74.4543\\
-1.7053 & 0.7777 & 13.9397 & -10.5308\\
-0.8277 & -0.3692 & 14.6466 & -52.9920\\
-0.8283 & -0.3694 & -0.6539 & -37.6904
\end{array}
\right]
\]
and $\mbox{spec}(\hat{A})=\{-40.0009,0.1620\pm0.8438i,15.3008\}$.

For another realization of noise with amplitude $10^{-2}$, we get the estimate
\end{example}

\[
\hat{A}=\left[
\begin{array}
[c]{cccc}%
-138.0893 & -60.7052 & -105.8111 & 301.5029\\
-0.2435 & 0.9101 & 12.9638 & -12.6745\\
-71.1408 & -31.9557 & -40.3842 & 142.3170\\
-71.1408 & -31.9557 & -55.6843 & 157.6172
\end{array}
\right]
\]
and $\mbox{spec}(\hat{A})=\{-40.1391,3.9326,0.9601,15.3002\}$.



The algorithm introduced above also allows us to compute the topological
entropy of linear systems, since it is determined by the unstable eigenvalues.
Recall that the topological entropy of a linear map on $\mathbb{R}^{n}$ is
defined in the following way:

Fix a compact subset $K\subset\mathbb{R}^{n}$, a time $\tau\in\mathbb{N}$ and
a constant $\varepsilon>0$. Then a set $R\subset\mathbb{R}^{n}$ is called
$(\tau,\varepsilon)$-spanning for $K$ if for every $y\in K$ there is $x\in R$
with%
\[
\left\Vert A^{j}y-A^{j}x\right\Vert <\varepsilon\text{ for all }j=0,...,\tau.
\]
By compactness of $K$, there are finite $(\tau,\varepsilon)$-spanning sets.
Let $R$ be a $(\tau,\varepsilon)$-spanning set of minimal cardinality
$\#R=r_{\min}(\tau,\varepsilon,K)$. Then%
\[
h_{top}(K,A,\varepsilon):=\lim_{\tau\rightarrow\infty}\frac{1}{\tau}\log
r_{\min}(\tau,\varepsilon,K),h_{top}(K,A):=\underset{\varepsilon
\rightarrow0^{+}}{\lim}h_{top}(K,\varepsilon).
\]
(the limits exist). Finally, the topological entropy of $A$ is%
\[
h_{top}(A):=\sup_{K}h_{top}(K,A),
\]
where the supremum is taken over all compact subsets $K$ of $\mathbb{R}^{n}.$

A classical result due to Bowen (cf. \cite[Theorem 8.14]{Walt82}) shows that
the topological entropy is determined by the sum of the unstable eigenvalues,
i.e.,%
\[
h_{top}(A)=\sum\max(1,\left\vert \lambda\right\vert ),
\]
where summation is over all eigenvalues of $A$ counted according to their
algebraic multiplicity.

Hence, when we approximate the unstable eigenvalues of $A$ by those of the
matrix $\hat{A}$, we also get an approximation of the topological entropy.

\begin{example}
\label{Example8}For Example \ref{Example6}, we get that $h_{top}(A)=34.80$
while for the estimate $\hat{A}$ one obtains $h_{top}(\hat{A})=34.7999$. For
Example \ref{Example7}, we get that $h_{top}(A)=55.30$ and $h_{top}(\hat
{A})=55.3008$. These estimates appear reasonably good.
\end{example}

\section{Identification of Linear Control Systems\label{Section4b}}

Consider the linear control system
\[
\label{lincontrol}x(k+1)=Ax(k)+Bu(k),
\]
with $A\in\mathbb{R}^{n\times n}$ and $B\in\mathbb{R}^{n\times1}$. We want to
estimate the matrices $A$ and $B$ from the time series $x(1)+\eta_{1},$
$\cdots,x(N)+\eta_{N}$ where $\eta$ satisfies the \textbf{Assumption} in
Section \ref{Section2}. The initial condition $x(0)$ and the control sequence
$u(0),$ $\cdots,u(N)$ are assumed to be known.

In order to estimate $A$ and $B$, we will extend algorithm $\mathcal{A}$. The
$i$th component of system (\ref{lincontrol}) is given by
\[
\label{lincontrol_i}x_{i}(k+1)=\sum_{j=1}^{n}a_{ij}x_{j}(k)+b_{i}u(k).
\]
For every $i$ we want to estimate the coefficients $b_{i}$ and $a_{ij},j=1,$
$\cdots,n$. Thus the linear map $f_{i}:\mathbb{R}^{n}\rightarrow\mathbb{R}$
given by
\[
\label{unknown_map1}(x_{1},...,x_{n},u)\mapsto\sum_{j=1}^{n}a_{ij}x_{j}%
+b_{i}u
\]
is unknown.
To extend algorithm $\mathcal{A}$, we will view system (\ref{lincontrol_i}) as
a system of the form (\ref{param1}) where the state $x$ is the extended state
$\underline{x}=(x,u)\in\mathbb{R}^{n}\times\mathbb{R}$ for (\ref{lincontrol}).
Hence, the kernel expansion (\ref{cij}) becomes%

\[
\label{cij1} {x}_{i}(k+1)=\sum_{j=1}^{N}c_{ij}\langle\underline{x}%
(j),\underline{x}(k)\rangle
\]
where $\underline{x}_{n+1}=u$ and the $c_{ij}$ satisfy the following set of equations:%

\begin{equation}
\label{seteqn1}\left[
\begin{array}
[c]{c}%
x_{i}(1)\\
\vdots\\
x_{i}(N)
\end{array}
\right]  =\Bigg(N\lambda I_{d}+\underline{\mathbb{K}}\Bigg)\left[
\begin{array}
[c]{c}%
c_{i1}\\
\vdots\\
c_{iN}%
\end{array}
\right]  ,
\end{equation}
with
\[
\underline{\mathbb{K}}=\left[
\begin{array}
[c]{ccc}%
\sum_{\ell=1}^{n+1}\underline{x}_{\ell}(1)\underline{x}_{\ell}(0) & \cdots &
\sum_{\ell=1}^{n+1}\underline{x}_{\ell}(N)\underline{x}_{\ell}(0)\\
\vdots & \cdots & \vdots\\
\sum_{\ell=1}^{n+1}\underline{x}_{\ell}(1)\underline{x}_{\ell}(N-1) & \cdots &
\sum_{\ell=1}^{n+1}\underline{x}_{\ell}(N)\underline{x}_{\ell}(N-1)
\end{array}
\right]  .
\]
Let us emphasize that $u=x_{n+1}$ does not appear on the left hand side of
(\ref{cij1})-(\ref{seteqn1}).

In reference to the case when $A$ has eigenvalues outside the unit circle, we
adopt the same method as in Section \ref{Section3} and define
\[
\label{gamma1}\underline{\sigma}:=\max\left\{  \frac{\left\Vert \underline
{x}(k+1)\right\Vert }{\left\Vert \underline{x}(k)\right\Vert },k\in
\{0,1,...,N\}\right\}  .
\]

\begin{example}
\label{Example9}(One Dimensional Case) Consider $x(k+1)=-0.9x(k)+3.5u$. For an
input $u(k)=\sin(k)+\cos(k)$ and for $100$ points we obtain the estimate
$\hat{A}=-0.9$ and $\hat{B}=3.5$ when there is no noise $\eta_{i}$. Here cross
validation gives $\gamma_{1}=1.5259\cdot10^{-05}$ and $\gamma_{2}=1$. For a
certain realization of the noise $\eta_{i}$ with amplitude $0.1$, we get
$\hat{A}=-0.9008$ and $\hat{B}=3.4983$. Here cross validation gives
$\gamma_{1}=0.0078$ and $\gamma_{2}=1$.
\end{example}

\begin{example}
\label{Example10}\ (Three Dimensional Stable Case) Consider control system
(\ref{lincontrol}) with
\[
A=\left[
\begin{array}
[c]{ccc}%
-0.9 & 1 & 0\\
0 & -0.1 & 1\\
0 & 0 & 0.8
\end{array}
\right]  \text{ and }B=\left[
\begin{array}
[c]{c}%
-2.5\\
-3.5\\
4.5
\end{array}
\right]  .
\]
With the input $u(k)=\sin(k)+\cos(k)$ and $100$ points, one computes the
estimates%
\[
\hat{A}=\left[
\begin{array}
[c]{ccc}%
-0.9000 & 1.0000 & 0.0000\\
0.0000 & -0.1000 & 1.0000\\
-0.0000 & -0.0000 & 0.8000
\end{array}
\right]  \text{ and }\hat{B}=\left[
\begin{array}
[c]{c}%
-2.5000\\
-3.5000\\
4.5000
\end{array}
\right]  .
\]
Here cross validation gives the regularization parameters $\gamma
_{i}=0.1526\cdot10^{-4}$ for $i=1,\cdots,4$. For some realization of
perturbations $\eta_{i}$ with amplitude $0.1$, one computes the estimates
\[
\hat{A}=\left[
\begin{array}
[c]{ccc}%
-0.9047 & 0.9984 & -0.0029\\
-0.0047 & -0.1016 & 0.9971\\
-0.0048 & -0.0018 & 0.7971
\end{array}
\right]  \text{ and }\hat{B}=\left[
\begin{array}
[c]{c}%
-2.5326\\
-3.5321\\
4.4661
\end{array}
\right]  .
\]
Here cross validation gives $\gamma_{1}=9.7656\cdot10^{-4}$, $\gamma
_{2}=9.7656\cdot10^{-4}$, $\gamma_{3}=1.5259\cdot10^{-5}$, $\gamma_{4}=4$.
\end{example}

\begin{example}
\label{Example11} (Three Dimensional Unstable Case) Consider control system
(\ref{lincontrol}) with
\[
A=\left[
\begin{array}
[c]{ccc}%
-20 & 1 & 0\\
0 & 1 & 1\\
0 & 0 & 20
\end{array}
\right]  \text{ and }B=\left[
\begin{array}
[c]{c}%
1\\
2\\
3
\end{array}
\right]  .
\]
The input $u(k)=\sin(k)+\cos(k)$ and $100$ points give the estimates%
\[
\hat{A}=\left[
\begin{array}
[c]{ccc}%
-19.9945 & 1.0009 & -0.0137\\
0.0013 & 0.9995 & 0.9919\\
0.0155 & -0.0171 & 19.7835
\end{array}
\right]  \text{ and }\hat{B}=\left[
\begin{array}
[c]{c}%
0.9898\\
1.9898\\
2.9333
\end{array}
\right]  .
\]
Here cross validation yields the regularization parameters $\gamma
_{i}=0.8882\cdot10^{-15}$ for $i=1,\cdots,4$. For some realization of
perturbations $\eta_{i}$ with amplitude $10^{-4}$, one computes the estimates
\[
\hat{A}=\left[
\begin{array}
[c]{ccc}%
-20.0000 & 0.9334 & -0.0058\\
-0.0008 & 0.9382 & 0.9939\\
-0.0008 & -0.0590 & 19.9937
\end{array}
\right]  \text{ and }\hat{B}=\left[
\begin{array}
[c]{c}%
0.9819\\
1.9814\\
2.9811
\end{array}
\right]  .
\]
Here cross validation gives $\gamma_{1}=\gamma_{2}=0.2384\cdot10^{-6}$,
$\gamma_{3}=\gamma_{4}=0.0596\cdot10^{-6}$.
\end{example}

These results show that algorithm $\mathcal{A}$ works quite well in these cases.

\section{Stabilization via Linear-Quadratic Optimal Control\label{Section5}}

A basic problem for linear control systems is stabilization by state feedback.
A standard method is to use linear quadratic optimal control, where the
feedback is computed using the solution of an algebraic Riccati equation. In
this section, we propose to replace in the algebraic Riccati equation the
system matrix $A$ by the estimate $\hat{A}$ obtained by learning theory.

The linear quadratic optimal control problem has the following form:\medskip

Minimize over all (continuous) inputs $u$%
\[
J_{\infty}(x_{0};u)=\sum_{k=0}^{\infty}\left[  x(k)^{\top}Qx(k)+u(k)^{\top
}Ru(t)\right]
\]
with $x(\cdot)$ given by
\[
{x}(k+1)=Ax(k)+Bu(k),\ k\geq0,\ x(0)=x_{0};
\]
here $Q\in\mathbb{R}^{n\times n}$ is positive semidefinite and $R\in
\mathbb{R}^{m\times m}$ is positive definite, and $A\in\mathbb{R}^{n\times
n},B\in\mathbb{R}^{n\times m}$.

Consider the discrete algebraic Riccati equation DARE
\[
A^{\top}(P-PB(R+B^{T}PB)^{-1}B^{\top}P)A+Q=P.
\]
Obviously, every solution $P$ is positive semi-definite. We cite the following
theorem from \cite{antsalkis}.

\textbf{Theorem. }Suppose that for every $x_{0}\in\mathbb{R}^{n}$ there is an
input $u$, such that $J(x_{0},u)<\infty$. Then the following holds:

(i) There is a unique solution $P$ of the DARE.

(ii) For every $x_{0}\in\mathbb{R}^{n}$ one has $J^{\ast}(x_{0}):=\inf
\{J(x_{0},u)\left\vert {}\right.  u$ an input$\}=x_{0}^{\top}Px_{0}$ and there
is a unique optimal input $u^{\ast}$ with $J^{\ast}(x_{0})=J(x_{0},u^{\ast})$.
This optimal input is generated by the feedback $F=(R+B^{T}PB)^{-1}B^{\top}PA$
and%
\[
u(k)=-Fx(k),k\geq0\text{.}%
\]
In particular, the feedback $F$ stabilizes the system, i.e., ${x}%
(k+1)=(A-BF)x(k)$ is stable.

Now we use an estimate $\hat{A}$ and $\hat{B}$ (obtained by kernel methods)
instead of $A$ and $B$ in the algebraic Riccati equation and obtain the
solution $\hat{P}$. Will the corresponding feedback $u=\hat{F}x:=-B^{\top}%
\hat{P}x$ also stabilize the system, i.e., is the following system stable:
\[
{x}(k+1)=(A-BB^{\top}\hat{P})x(k)?
\]

\begin{example}
\label{Example12}Consider the one-dimensional system $x(k+1)=-0.9x(k)+3.5u$ in
Example \ref{Example9}. In the absence of noise, we get $\hat{A}=-0.9$ and
$\hat{B}=3.5$. We have that $A-B\hat{F}=\hat{A}-\hat{B}\hat{F}=-0.0643$. When
there is noise of amplitude $0.1$, we get that $\hat{A}=-0.9002$ and $\hat
{B}=3.4929$ and $A-B\hat{F}=-0.0643$ while $\hat{A}-\hat{B}\hat{F}=-0.0610$.
Hence, the controller improves stability.
\end{example}

\begin{example}
\label{Example13}Consider control system (\ref{lincontrol}) with%
\[
A=\left[
\begin{array}
[c]{ccc}%
-0.9 & 1 & 0\\
0 & -0.1 & 1\\
0 & 0 & 0.8
\end{array}
\right]  \text{ and }B=\left[
\begin{array}
[c]{c}%
-2.5\\
-3.5\\
4.5
\end{array}
\right]  .
\]
As illustrated in Example \ref{Example10}, without noise we get excellent
approximations of $A$ and $B$. For both cases, the set of eigenvalues of the
closed-loop system is $\{-0.6172,0.4049,-0.0018\}$. With a noise of maximal
amplitude $0.1$, the estimates $\hat{A}$ and $\hat{B}$ are given in Example
\ref{Example10}. For the feedback system one finds
\begin{align*}
\mbox{spec}(\hat{A}-\hat{B}\hat{F})  &  =\{-0.6204,0.4053,-0.0018\},\\
\mbox{spec}(A-{B}\hat{F})  &  =\{-0.6240,-0.0062,0.4111\}.
\end{align*}
In this example the feedback based on the estimate also stabilizes the
original system.
\end{example}

\begin{example}
\label{Example14} Consider control system (\ref{lincontrol}) with%
\[
A=\left[
\begin{array}
[c]{ccc}%
-20 & 1 & 0\\
0 & 1 & 1\\
0 & 0 & 20
\end{array}
\right]  \text{ and }B=\left[
\begin{array}
[c]{c}%
1\\
2\\
3
\end{array}
\right]  .
\]
As Example \ref{Example11} illustrates, without noise we get excellent
approximations of $A$ and $B$. For the feedback system one finds%
\begin{align*}
\mbox{spec}(\hat{A}-\hat{B}\hat{F})  &  =\{0.1994,0.0483,-0.0501\},\\
\mbox{spec}(A-{B}\hat{F})  &  =\{-0.1234\pm2.0777i,0.5279\}.
\end{align*}
When there is noise of amplitude $10^{-4}$, one computes the estimates%
\[
\hat{A}=\left[
\begin{array}
[c]{ccc}%
-19.9805 & 0.7484 & 0.0135\\
-0.0062 & 0.7969 & 1.0107\\
-0.0229 & 0.9851 & 19.6776
\end{array}
\right]  \text{ and }\hat{B}=\left[
\begin{array}
[c]{c}%
1.0194\\
2.0114\\
2.6673
\end{array}
\right]  .
\]
This are bad approximations for $A$ and $B$. Furthermore, for the feedback
system one finds
\begin{align*}
\mbox{spec}(\hat{A}-\hat{B}\hat{F})  &  =\{0.1929,0.0477,-0.0501\},\\
\mbox{spec}(A-{B}\hat{F})  &  =\{1.4510\pm3.0103i,-2.5232\}.
\end{align*}
Thus the stabilizing controller for the approximate system does not stabilize
the true system.
\end{example}

\section{Conclusions\label{Conclusions}}

This paper has introduced the algorithm $\mathcal{A}$ based on kernel methods
to identify a stable linear dynamical system from a time series. The numerical
experiments give excellent results in the absence of noise and structural
perturbations. In the presence of noise and structural perturbations the
algorithm works well in the stable case. In the unstable case, a modified
algorithm works quite well in the presence of noise but cannot handle
structural perturbations.

Then we have extended algorithm $\mathcal{A}$ to identify linear control
systems. In particular, we have used estimates obtained by kernel methods to
stabilize linear systems using linear-quadratic control and the algebraic
Riccati equation. Here the numerical experiments seem to indicate that the
same conclusions on applicability of the algorithm apply.

Extensions of the considered algorithms to nonlinear systems appear feasible
and are left to future work.

\appendix

\section{Appendix: Elements of Learning Theory}

In this section, we give a brief overview of Reproducing Kernel Hilbert Spaces
(RKHS) as used in statistical learning theory. The discussion here borrows
heavily from~Cucker and Smale \cite{cucker}, Wahba \cite{Wahba}, and
Sch\"{o}lkopf and Smola \cite{smola_book}. Early work developing the theory of
RKHS was undertaken by I.J. Schoenberg \cite{schoenberg, schoenberg1,
schoenberg2} and then N. Aronszajn~\cite{AronRKHS}. Historically, RKHS came
from the question, when it is possible to embed a metric space into a Hilbert space.

\begin{definition}
Let $\mathcal{H}$ be a Hilbert space of functions on a set $\mathcal{X}$ which
is a closed subset of $\mathbb{R}^{n}$. Denote by $\langle f, g \rangle$ the
inner product on $\mathcal{H}$ and let $||f||= \langle f, f \rangle^{1/2}$ be
the norm in $\mathcal{H}$, for $f$ and $g \in\mathcal{H}$. We say that
$\mathcal{H}$ is a reproducing kernel Hilbert space (RKHS) if there exists
$K:\mathcal{X} \times\mathcal{X} \rightarrow\mathbb{R}$ such that

\begin{itemize}
\item[i.] $K$ has the reproducing property, i.e., $f(x)=\langle f(\cdot
),K(\cdot,x)\rangle$ for all $f\in\mathcal{H}$.

\item[ii.] $K$ spans $\mathcal{H}$, i.e., $\mathcal{H}=\overline
{\mbox{span}\{K(x,\cdot)|x\in\mathcal{X}\}}$.
\end{itemize}

$K$ will be called a reproducing kernel of $\mathcal{H}$ and $\mathcal{H}_{K}$
will denote the RKHS $\mathcal{H}$ with reproducing kernel $K$.
\end{definition}

\begin{definition}
Given a kernel $K:\mathcal{X}\times\mathcal{X}\rightarrow\mathbb{R}$ and
inputs $x_{1},\cdots,x_{n}\in\mathcal{X}$, the $n\times n$ matrix
\[
k:=(K(x_{i},x_{j}))_{ij},
\]
is called the \emph{Gram Matrix} of $k$ with respect to $x_{1},\cdots,x_{n}$.
If for all $n\in\mathbb{N}$ and distinct $x_{i}\in\mathcal{X}$ the kernel
$K\ $gives rise to a strictly positive definite Gram matrix, it is called
strictly positive definite.
\end{definition}

\begin{definition}
(Mercer kernel map) A function $K:\mathcal{X} \times\mathcal{X} \rightarrow
\mathbb{R}$ is called a Mercer kernel if it is continuous, symmetric and
positive definite.
\end{definition}

The important properties of reproducing kernels are summarized in the
following proposition.

\begin{proposition}
\label{prop1} If $K$ is a reproducing kernel of a Hilbert space $\mathcal{H}$, then

\begin{itemize}
\item[i.] $K(x,y)$ is unique.

\item[ii.] For all $x,y\in\mathcal{X}$, $K(x,y)=K(y,x)$ (symmetry).

\item[iii.] $\sum_{i,j=1}^{m}\alpha_{i}\alpha_{j}K(x_{i},x_{j}) \ge0$ for
$\alpha_{i} \in\mathbb{R}$ and $x_{i} \in\mathcal{X} $ (positive definitness).

\item[iv.] $\langle K(x,\cdot),K(y,\cdot) \rangle_{\mathcal{H}}=K(x,y)$.

\item[v.] The following kernels, defined on a compact domain $\mathcal{X}%
\subset\mathbb{R}^{n}$, are Mercer kernels: $K(x,y)=x\cdot y^{\top}$ (Linear),
$K(x,y)=(1+x\cdot y^{\top})^{d},\quad d\in\mathbb{N}$ (Polynomial),
$K(x,y)=e^{-\frac{||x-y||^{2}}{\sigma^{2}}},\quad\sigma>0$ (Gaussian).
\end{itemize}
\end{proposition}

\begin{theorem}
\label{thm1} Let $K:\mathcal{X}\times\mathcal{X}\rightarrow\mathbb{R}$ be a
symmetric and positive definite function. Then there exists a Hilbert space of
functions $\mathcal{H}$ defined on $\mathcal{X}$ admitting $K$ as a
reproducing Kernel. Moreover, there exists a function $\Phi:X\rightarrow
\mathcal{H}$ such that
\[
K(x,y)=\langle\Phi(x),\Phi(y)\rangle_{\mathcal{H}}\quad\mbox{for}\quad
x,y\in\mathcal{X}.
\]
$\Phi$ is called a feature map.

Conversely, let $\mathcal{H}$ be a Hilbert space of functions $f:\mathcal{X}%
\rightarrow\mathbb{R}$, with $\mathcal{X}$ compact, satisfying
\[
\text{For all }x\in\mathcal{X}\text{ there is }\kappa_{x}%
>0,~\mbox{such that}~|f(x)|\leq\kappa_{x}||f||_{\mathcal{H}}.
\]
Then $\mathcal{H}$ has a reproducing kernel $K$.
\end{theorem}

\begin{remarks}

\begin{itemize}
\item[i.] The dimension of the RKHS can be infinite and corresponds to the
dimension of the eigenspace of the integral operator $L_{K}:\mathcal{L}_{\nu
}^{2}(\mathcal{X})\rightarrow\mathcal{C}(\mathcal{X})$ defined as
$(L_{K}f)(x)=\int K(x,t)f(t)d\nu(t)$ if $K$ is a Mercer kernel, for
$f\in\mathcal{L}_{\nu}^{2}(\mathcal{X})$ and $\nu$ a Borel measure on
$\mathcal{X}$.

\item[ii.] In Theorem \ref{thm1}, and using property [iv.] in Proposition
\ref{prop1}, we can take $\Phi(x):=K_{x}:=K(x,\cdot)$ in which case
$\mathcal{F}=\mathcal{H}$ -- the \textquotedblleft feature
space\textquotedblright\ is the RKHS. This is called the \emph{canonical
feature map}.

\item[iii.] The fact that Mercer kernels are positive definite and symmetric
shows that kernels can be viewed as generalized Gramians and covariance matrices.

\item[iv.] In practice, we choose a Mercer kernel, such as the ones in [v.] in
Proposition \ref{prop1}, and Theorem \ref{thm1}, that guarantees the existence
of a Hilbert space admitting such a function as a reproducing kernel.
\end{itemize}
\end{remarks}


RKHS play an important role in learning theory whose objective is to find an
unknown function
\[
\label{unknown_func} f^{\ast}:X\rightarrow Y
\]
from random samples%

\[
\label{samples1} \mathbf{s}=(x_{i},y_{i})|_{i=1}^{m},
\]

In the following we review results from \cite{smale_shannon1} (for a more
general setting, cf. \cite{cucker}) in the special case when the data samples
$\mathbf{s}$ are such that the following assumption holds.

\emph{Assumption 1:} The samples in (\ref{samples1}) have the special form
\[
\mathcal{S: \quad}\label{samples} \mathbf{s}=(x,y_{x})|_{x \in\bar{x}},
\]
where $\bar{x}=\{x_{i}\}|_{i=1}^{d+1}$ and $y_{x}$ is drawn at random from
$f^{\ast}(x)+\eta_{x} $, where $\eta_{x}$ is drawn from a probability measure
$\rho_{x}$.

Here for each $x \in X$, $\rho_{x}$ is a probability measure with zero mean,
and its variance $\sigma_{x}^{2}$ satisfies $\sigma^{2} :=\sum_{x \in\bar{x}}
\sigma_{x}^{2} < \infty$. Let $X$ be a closed subset of $\mathbb{R}^{n}$ and
$\bar{t} \subset X$ is a discrete subset. Now, consider a kernel $K: X \times
X \rightarrow\mathbb{R}$ and define a matrix (possibly infinite) $K_{\bar
{t},\bar{t}} : \ell^{2}(\bar{t}) \rightarrow\ell^{2}(\bar{t})$ as
\[
(K_{\bar{t},\bar{t}}a)_{s} = \sum_{t \in\bar{t}}K(s,t)a_{t}, \quad s \in
\bar{t}, a \in\ell^{2}(\bar{t}),
\]
where $\ell^{2}(\bar{t})$ is the set of sequences $a=(a_{t})_{t \in\bar{t}}:
\bar{t} \rightarrow\mathbb{R}$ with $\langle a,b \rangle=\sum_{t \in\bar{t}%
}a_{t} b_{t}$ defining an inner product. For example, we can take
$X=\mathbb{R}$ and $\bar{t}=\{0,1,\cdots,d\}$.

In the case of a linear dynamical system (\ref{syst1}), we are interested in
learning the map $x(k)\mapsto x(k+1)$. Here we can apply the following results.

The problem to approximate a function $f^{\ast}\in\mathcal{H}_{K}$ from
samples $\mathbf{s}$ of the form (\ref{samples1}) has been studied in
\cite{smale_shannon1, smale_shannon2}. It is reformulated as the minimization
problem
\[
\label{reg_opt1}\bar{f}_{\mathbf{s},\gamma}:=\mbox{arg}{\mbox{min}}%
_{f\in\mathcal{H}_{K,\bar{t}}}\bigg\{\sum_{x\in\bar{x}}(f(x)-y_{x})^{2}%
+\gamma||f||_{K}^{2}\bigg\},
\]
where $\gamma\geq0$ is a regularization parameter. Moreover,when $\bar{x}$ is
not defined by a uniform grid on $X$, the authors of \cite{smale_shannon1}
introduced a weighting $w:=\{w_{x}\}_{x\in\bar{x}}$ on $\bar{x}$ with
$w_{x}>0$\footnote{A suggestion in \cite{smale_shannon1} is to consider the
$\rho_{X}-$volume of the Voronoi cell associated with $\bar{x}$. Another
example is $w=1$ or if $|\bar{x}|=m<\infty$, $w=\frac{1}{m}$.}. Let $D_{w}$ be
the diagonal matrix with diagonal entries $\{w_{x}\}_{x\in\bar{x}}$. Then,
$||D_{w}||\leq||w||_{\infty}$.

In this case, the regularization scheme (\ref{reg_opt1}) becomes
\[
\label{reg_opt2}\bar{f}_{\mathbf{s},\gamma}:=\mbox{arg}{\mbox{min}}%
_{f\in\mathcal{H}_{K,\bar{t}}}\bigg\{\sum_{x\in\bar{x}}w_{x}(f(x)-y_{x}%
)^{2}+\gamma||f||_{K}^{2}\bigg\},
\]

\begin{theorem}
\label{thm_sol1} Assume $f^{\ast}\in\mathcal{H}_{K,\bar{t}}$ and the standing
hypotheses with $X$, $K$, $\bar{t}$, $\rho$ as above, $y$ as in (\ref{samples}%
). Suppose $K_{\bar{t},\bar{x}}D_{w}K_{\bar{x},\bar{t}}+\gamma K_{\bar{t}%
,\bar{t}}$ is invertible. Define $\mathcal{L}$ to be the linear operator
$\mathcal{L}=(K_{\bar{t},\bar{x}}D_{w}K_{\bar{x},\bar{t}}+\gamma K_{\bar
{t},\bar{t}})^{-1}K_{\bar{t},\bar{x}}D_{w}$. Then problem (\ref{reg_opt2}) has
the unique solution
\[
f_{\mathbf{s},\gamma}=\sum_{t\in\bar{t}}(\mathcal{L}y)_{t}K_{t}%
\]

\end{theorem}

\emph{Assumption 2}: For each $x \in X$, $\rho_{x}$ is a probability measure
with zero mean supported on $[-M_{x},M_{x}]$ with $\mathcal{B}_{w} :=(\sum_{x
\in\bar{x}}w_{x} M_{x}^{2})^{\frac{1}{2}} < \infty$.

The next theorems give estimates for the different sources of errors.

\begin{theorem}
\label{TheoremA.3}(Sample Error) \cite[Theorem 4, Propositions 2 and
3]{smale_shannon1} Let \emph{Assumptions 1} and \emph{2} be satisfied, suppose
that $K_{\bar{t},\bar{x}}D_{w}K_{\bar{x},\bar{t}}+\gamma K_{\bar{t},\bar{t}}$
is invertible and let $f_{\mathbf{s},\gamma}=\sum_{t\in\bar{t}}c_{t}K_{t}$ be
the solution of (\ref{reg_opt2}) given in Theorem \ref{thm_sol1} by
$c=\mathcal{L}y$. Define
\begin{align*}
\mathcal{L}_{w}  &  :=(K_{\bar{t},\bar{x}}D_{w}K_{\bar{x},\bar{t}}+\gamma
K_{\bar{t},\bar{t}})^{-1}K_{\bar{t},\bar{x}}D_{w}^{1/2}\\
\kappa &  :=||K_{\bar{t},\bar{t}}||\;||(K_{\bar{t},\bar{x}}D_{w}K_{\bar
{x},\bar{t}}+\gamma K_{\bar{t},\bar{t}})^{-1}||^{2}.
\end{align*}
Then for every $0<\delta<1$, with probability at least $1-\delta$ we have the
sample error estimate
\[
\label{epsilon_samp}||f_{\mathbf{s},\gamma}-f_{\bar{x},\gamma}||_{K}^{2}%
\leq\mathcal{E}_{\mbox{samp}}:=\kappa\sigma_{w}^{2}\alpha^{-1}\bigg(\frac
{2||K_{\bar{t},\bar{t}}\mathcal{L}_{w}||\;||\mathcal{L}_{w}||\;\mathcal{B}%
_{w}^{2}}{\kappa\sigma_{w}^{2}}\;\log{\frac{1}{\delta}}\bigg),
\]
where $\alpha(u):=(u-1)\log u$ for $u>1$. In particular, $\mathcal{E}%
_{\mbox{samp}}\rightarrow0$ when $\gamma\rightarrow\infty$ or $\sigma_{w}%
^{2}\rightarrow0$.
\end{theorem}

\begin{theorem}
\label{TheoremA.4}(Regularization Error and Integration Error)
\cite[Proposition 4 and Theorem 5]{smale_shannon1} Let \emph{Assumptions 1}
and \emph{2} be satisfied and let $\bar{X}=(X_{x})_{x\in\bar{x}}$ be the
Voronoi cell of $X$ associated with $\bar{x}$ and $w_{x}=\rho_{X}(X_{x})$.
Define the Lipschitz norm on a subset $X^{\prime}\subset X$ as
$||f||_{\mbox{Lip}(X^{\prime})}:=||f||_{L^{\infty}(X^{\prime})}+\sup_{s,u\in
X}\frac{|f(s)-f(u)|}{||s-u||_{\ell^{\infty}(\mathbb{R}^{n})}}$ and assume that
the inclusion map of $\mathcal{H}_{K,\bar{t}}$ into the Lipschitz space
satisfies\footnote{This assumption is true if $X$ is compact and the inclusion
map of $\mathcal{H}_{K,\bar{t}}$ into the space of Lipschitz functions on $X$
is bounded which is the case when $K$ is a $C^{2}$ Mercer kernel
\cite{zhou_capacity}. In fact, if $||f||_{\mbox{Lip}(X)}\leq C_{0}||f||_{K}$
for each $f\in\mathcal{H}_{K,\bar{t}}$, then $C_{\bar{x}}\leq C_{0}^{2}%
\rho_{X}(X)$.}%

\[
C_{\bar{x}}:=\sup_{f\in\mathcal{H}_{K,\bar{t}}}\frac{\sum_{x\in\bar{x}}%
w_{x}||f||_{\mbox{Lip}(X_{x})}^{2}}{||f||_{K}^{2}}<\infty.
\]
Suppose that $\bar{x}$ is $\Delta-$dense in $X$, i.e., for each $y\in X$ there
is some $x\in\bar{x}$ satisfying $||x-y||_{\ell^{\infty}(\mathbb{R}^{n})}%
\leq\Delta$.

Then for $f^{\ast}\in\mathcal{H}_{K,\bar{t}}$
\[
||f_{\bar{x},\gamma}-f^{\ast}||^{2}\leq||f^{\ast}||_{K}^{2}(\gamma+8C_{\bar
{x}}\Delta)
\]

\end{theorem}

\begin{theorem}
(Sample, Regularization and Integration Errors) \label{thm:errors}
\cite[Corollary 5]{smale_shannon1} Under the assumptions of Theorems
\ref{TheoremA.3} and \ref{TheoremA.4}, let $\bar{X}=(X_{x})_{x\in\bar{x}}$ be
the Voronoi cell of $X$ associated with $\bar{x}$ and $w_{x}=\rho_{x}(X_{x}%
)$.
Suppose that $\bar{x}$ is $\Delta-$dense, $C_{\bar{x}}<\infty$, and $f^{\ast
}\in\mathcal{H}_{K,\bar{t}}$. Then, for every $0<\delta<1$, with probability
at least $1-\delta$ there holds
\[
\label{err_est1}||f_{\mathbf{{s},\gamma}}-f^{\ast}||^{2}\leq2C_{\bar{x}%
}\mathcal{E}_{\mbox{samp}}+2||f^{\ast}||_{K}^{2}(\gamma+8C_{\bar{x}}\Delta),
\]
where $\mathcal{E}_{\mbox{samp}}$ is given in (\ref{epsilon_samp}).
\end{theorem}

\section{Acknowledgements}

BH thanks the European Commission and the Scientific and the Technological
Research Council of Turkey (Tubitak) for financial support received through a
Marie Curie Fellowship.

\end{document}